\documentclass[12pt]{article}
\usepackage{amssymb}
\usepackage{graphics}
\newcommand\Af{\Bbb A}                        
\newcommand\fl{ f_{\lambda}}
\newcommand\e{\vspace{15pt}}                   

\renewcommand\Pr{ \Bbb P}                         
\newcommand\Z{\Bbb Z}                          
\newcommand\R{\Bbb R}                          
 
\newcommand\fin{\hfill $\diamond $}             
\newtheorem{Def}{\bf Definition}
\newtheorem{Obs}[Def]{\bf Remark}
\newtheorem{Prop}[Def]{\bf Proposition}
\newtheorem{Lemm}[Def]{\bf Lemma}
\newtheorem{Cor}[Def]{\bf Corollary}

\begin{document}
\title{Real Quartic  Surfaces in $\R{\Pr}_3 $ containing 16  Skew Lines}
\author{Isidro  Nieto \\ {\it Instituto de Fisica y Matematicas de la Universidad Michoacana}}

\date{}
\maketitle
\begin{abstract}
In [BN] the authors construct  a special complex of degree 20 over $M$, which for an open 
three dimensional set parametrizes smooth complex surfaces of degree four invariant which are Heisenberg invariant, and~each member of the family  contains 32 lines but only 16 skew lines. The coordinates of the lines however need not be real.~For a  point $l$  in a Zariski open set  of $ M $ an algorithm 
is  presented which  evaluates  the {\it  real } coefficients of  $l$  in terms of the K- coordinates of  $l$.~The author uses a code in Maple which allows  him to  construct  very explicit examples  of  real smooth Heisenberg invariant Kummer surfaces  containing  the special configuration  of 16 {\it real} skew lines.~An example is presented at the end of the paper. 
\end{abstract}
\section{ Introduction}\label{intro}
In the variables $ x_0, x_1, x_2 , x_3, x_4, x_5 $ consider
 the equation :
\begin{eqnarray}
x_0^2 + x_1^2 + x_2^2 + x_3^2 + x_4^2 + x_5^2  & = 0  \label{profe} \\
1/ x_0^2 + 1/ x_1^2 +  1/x_2^2 + 1/x_3^2 + 1/x_4^2 + 1/x_5^2 &  = 0. \label{proble}
\end{eqnarray}
This is  a rational form of the equation which must be written in regular form ! \\
 Let us consider a solution to these set of equations.
 For this let $u_i = x_i^2$ and  fixing $u_2, u_3, u_4, u_5$ then in \cite{Nie}
 we have found an equation of the form $ u_0^2 + \lambda u_0 + \lambda / \mu = 0 $ where $ \lambda, \mu$ are
 constants depending on the above fixed variables. The solutions  can however be {\it complex}.~If we let $\iota = \sqrt{-1} $ and 
introducing
 $$
\begin{array}{cccccc}
p_{01} & = x_0 - \iota x_1,   & \quad p_{02} & = x_2 - \iota x_3 ,  & \quad p_{03} & = x_4 - \iota x_5 ,  \\
p_{23} & =- ( x_0 + \iota x_1), & \quad  p_{13} & = x_2 + \iota x_3 , & \quad  p_{12} & = -( x_4 + \iota x_5)
\end{array}
$$
  if all $\{ p_{i,j} \}$ are {\it real}, the $p_{i,j}$ coordinates will
 determine a real line $l$ in ${\Pr}_3$ contained in a real quartic surface, which
 is invariant under the well known group $H^t $ the{\it Heisenberg group} of level 2.~Such surfaces form a three dimensional parameter space within the thirty four
 dimensional space of quartic surfaces in $\R {\Pr}_3$ and are therefore very special .~Here we
present two  algorithms for writing the equation for the surface and the lines  lying on them together with
 one example illustrating the second algorithm.~These examples  are   surfaces   of degree four
  embedded as  Kummer surfaces associated to    abelian surfaces of polarization  type $(1,3)$ (c.f.~\cite{BN} and \cite{Nie}).
 We give an elementary  solution
 to the problem stated  in the beginning of  section \ref{algor} which solution is    the first algorithm.~In  section \ref{algorn} we introduce the second algorithm .~In section  \ref{superf} we  briefly recall  the basic  concepts of Line Geometry and
 Quartic Surfaces in $\R{\Pr}_3$ to solve our problem ,use the second algorithm 
  in section \ref{ejemp} to construct  one   example and  briefly sketch  how  the code was  written in Maple. Additionally  we produce the equation of the 
quartic surface and the set of skew lines lying on it.   
 \section{An elementary Problem} \label{algor}
We consider the following elementary problem which first solution was communicated by the author by W. Barth in an unpublished manuscript:
\noindent
{\it Find solutions to the set of  equations \ref{profe} and  \ref{proble}   such that  $ x_0, x_2, x_4 $ are real
 and $ x_1, x_3, x_5 $ are purely imaginary.}

\noindent For later use we  let
 \begin{eqnarray}
x_1 = \iota y_1 , x_3 = \iota y_3 ,  x_5 = \iota y_5 \label{problema} 
 \end{eqnarray}
 where $\iota$ is $ \sqrt{-1} $.
 By substitution the new form of equation~\ref{profe} is
\begin{eqnarray}
x_0^2  + x_2^2 + x_4^2 =  y_1^2 + y_3^2 + y_5^2  \label{equad}
\end{eqnarray}
and similarly for the second equation
 $$
1/x_0^2 + 1/x_2^2 + 1/x_4^2 = 1/y_1^2 + 1/y_3^2 + 1/y_5^2 .
$$
We will normalize  equation ~\ref{equad}  as
$$
x_0^2 + x_2^2 + x_4^2 = 1 = y_1^2 + y_3^2 + y_5^2
$$
and set
$$
q= 1/ x_0^2 + 1/x_4^2 + 1/x_2^2 = 1/ y_1^2 + 1/y_3^2 + 1/ y_5^2 .
$$
If we introduce the following change of variables:
$$
q_i = \left \{ \begin{array}{cc}
 x_i^2  & \mbox{for} \, \, \mbox{i  even}, \\
 y_i^2  & \mbox{otherwise.}
\end{array}
\right.
$$
we obtain the following
\begin{Prop}
 {\it Choose }  $ q_4 \neq q_5 $ {\it real such that }
$$
0 \lneq  q_4, q_5 \lneq  1
$$
 and
 $$
 q \gneq  \mbox{max} ( 1/q_4 + 4/( 1-q_4), 1/ q_5 + 4/(1-q_5))
 $$
{\it then there exist}
$ x_0, x_2, x_4, y_1, y_3, y_5 $ {\it real satisfying  equation~\ref{proble}}.
\end{Prop}
{\it Proof}.-~Considering
$$
1/q_0 + 1/ q_2 = q - 1/q_4
$$
and setting $ s=  1 - q_4= q_0 + q_2 $ we obtain the formal equation
 $$
\begin{array}{cc}
 p = q_0 q_2 & =   ( q_0 + q_2 ) / ( q - 1/q_4), \\ \nonumber
             &  =   q_4 s/ ( q q_4 - 1), \\ \nonumber 
             & = q_4(1-q_4)/( q q_4 -1).     \label{equac}
\end{array}
$$
\noindent Assuming that $ s, p $ are given then the roots of
\begin{eqnarray}
X^2 - s X + p = 0  \label{equat}
\end{eqnarray}
determine  $ q_0, q_2 $ whenever  $ q$ is given. In fact, $ q_0 = ( s + \sqrt{ s^2 - 4 p})/2 , \linebreak
 q_2=  ( s-\sqrt{s^2 -4p})/2 $ . 
\e
\noindent Similarly  if $ s' = 1 - q_5= q_1 + q_3 $ we obtain  the quadratic equation:
$$
X^2 - s' X + p' = 0
$$
and the condition  $ q_4 \neq q_5 $ implies that  $s \neq s' $.~Equation~\ref{equat}
 has  the elementary  solutions
$$
( s \pm \sqrt{ s^2 - 4 p } )/2 .
$$
\noindent  A necessary condition  that  equation \ref{equat} has different solutions is that
$$
s^2 \gneq 4 p = 4 s q_4/( q q_4 - 1).
$$
\noindent By assumption $ s \ge 0 $ thus  $ s \gneq 4 q_4/( q q_4 - 1) $ which gives
$$
q \ge 1/q_4 + 4/s.
 $$
\noindent It also follows from the last  inequality and  equation~\ref{equac}  that
$$
 p  \ge 0 \, \, \, q_0, q_2 \ge 0.
 $$
 \noindent Similarly if $ q \gneq 1/q_5 + 4/( 1-q_5) $ and  $ 0 \lneq q_5 \lneq 1 $ then  $ y_1, y_3 $ are roots of positive numbers.
\fin
\begin{Obs} For this condition the normalization
 of   equation~\ref{equad} is not needed. Moreover, the condition chosen for convenience $q_4= q_5$ is indeed necessary since it  implies that
$$
x_0^2 + x_1^2 = x_4^2 +x_5^2 = x_2^2 +x_3^2 =0
$$
and it follows from [\cite{Nie1}, \S 5 and \S 6] that the surface is elliptic ruled or singular along four nodes.
\end{Obs}
\section{The Algorithm}\label{algorn}
  Once again we fix the normalization of equation ~\ref{equad}  in the 
 coordinates $q_i $, which is  $ 1 = q_0 + q_2 + q_4 $. We are to find real positive  solutions to  the non-linear  system of polynomial equations given in rational form as :
\begin{eqnarray}
 \lambda & =&1/q_2 + 1/q_0 + 1/q_4  \nonumber \\
  1    & = &  q_0 + q_2 + q_4 \label{ecuaci}
\end{eqnarray}
\noindent For the following we let $ \omega = (\lambda q_0^2 -q_0( \lambda -3)  +1)/\lambda $ and 
$N = ( \lambda q_0 -1)(1-q_0) $, $ \Delta^2 = N ( N -4q_0) $.~Note that   these are polynomial expressions in $ q_0 $. We fix in the sequel  the numerical value for $q_0 $.
\begin{Prop}\label{qqq}
The roots  of equation~ \ref{ecuaci} for $ q_2 $  are:
$ q_2=  (N  \pm \Delta) / 2(\lambda q_0 -1)$ with  $\Delta^2 = \lambda^2 \omega ( 1/\lambda - q_0)(1-q_0) $. 
\end{Prop}
\noindent{\it Proof}.- Claim: 
$$
q_2 = ( \lambda q_0 -1)(1-q_0) \pm \sqrt{ (\lambda q_0 -1)^2 (1-q_0)^2 - 4 q_0(1-q_0)(\lambda q_0 -1)} / 2( \lambda q_0 -1)
$$
\noindent {\it To prove the claim}  from the above given  equation one obtains 
$$
 q_0 q_4 + q_2q_4 + q_0q_2 - \lambda q_0 q_2 q_4 = 1- q_0 -q_2 -q_4 =0 
$$
\noindent By substituting  the value of $q_4$ in the last equation  one obtains the quadratic equation 
$$
q_2^2( \lambda q_0 -1) -q_2( \lambda q_0 -1) ( 1-q_0)  +  q_0 ( 1-q_0) =0 
$$
\noindent which solution  is as claimed.

\noindent To finish the proof  we write $ \Delta^2 = N M $ where $M= ( \lambda q_0 -1)(1-q_0) -4 q_0 $.~The last 
equality is equal to $ -( \lambda q_0^2 +  q_0( \lambda -3) + 1) $. 
\fin 
\begin{Obs}
 The	condition $ 1= q_0 + q_2
 + q_4 $ still has to be fulfilled which does not 
follow from Proposition \ref{qqq}.
\end{Obs}
\begin{Obs}\label{rema}
 Let $ q_0 >  1/\lambda $  then  
$$
q_0 + q_2 < 1 \Longleftrightarrow N  - \epsilon \Delta > 0 \,\, \mbox{with}
 \,\, \epsilon \in \{ \pm 1 \} .
$$
\end{Obs}
\noindent{\it Proof of the remark}~An easy computation shows that 
$$ 
q_0 + q_2 = ((\lambda q_0 -1)(q_0 +
1) + \epsilon  \Delta) / 2( \lambda q_0 -1). 
$$
\noindent Hence obtain $ (\lambda q_0 -1)(2-( q_0 +1)) +  \epsilon  \Delta >0 $.
\fin  

\noindent  We introduce more notation. Let $ \sigma = 
 ( \lambda -3)/ 2 \lambda, \rho 
= (\lambda -9)(\lambda -1)/4\lambda^2 $.  

\begin{Prop}\label{imp} $ \Delta^2 > 0 \Longleftrightarrow   \omega >
 0 , q_0 < 1/\lambda   \,\,\,\, \mbox{or}\,\,\,\,
\omega < 0 ,
 1/\lambda < q_0$. 
\end{Prop}
\noindent{\it Proof}.-~This is a consequence of Proposition~\ref{qqq}.
\fin
\begin{Lemm}\label{lemu}$  \omega > 0 \Longleftrightarrow  q_0  > \sigma  + \sqrt{\rho} $ or 
$ q_0 < \sigma - \sqrt{\rho} $ .
\end{Lemm}
\noindent{\it Proof}.-~By writing:
$$
\begin{array}{cc} 
\omega   & =   q_0^2 -  q_0(\lambda -3)/\lambda + 1/\lambda,  \\
     & =   ( q_0^2 - 
2 (\lambda -3)/2\lambda q_0 + (\lambda -3)^2/ 4\lambda^2 ) + 1/\lambda - (\lambda -3)^2/4 \lambda^2, \\    
    &  =  ( q_0 -(\lambda -3)/2\lambda)^2 - \rho. 
\end{array}
$$
\noindent The cases claimed in the lemma follow inmediately from it.  
\fin

\noindent The proof of the following lemmma is analogous to the previous one:
\begin{Lemm}\label{valor} $ \omega < 0 \Longleftrightarrow  q_0 \in ( \sigma -\sqrt{\rho},   \sigma +
\sqrt{\rho} ) $.
\end{Lemm} 
\fin
\begin{Lemm}\label{lemp} If $N > 0 $ then $ N - \Delta > 0 $. If additionally  $ q_0 \in ( \sigma - \sqrt{\rho}, \sigma + \sqrt{\rho}) $ then $ q_0 + q_2<1 $. 
\end{Lemm}
\noindent{\it Proof}.- $\Delta^2  = - N \lambda  \omega $. We have to evaluate $  N - \sqrt{-N \lambda \omega} $. But  
$$
(N- \sqrt{ -N \lambda \omega } )(N + \sqrt{-N \lambda \omega} ) = N ( N   +\lambda \omega )  
$$
\noindent From  $ N + \lambda \omega  =  4 q_0 > 0 $ thus  $ N -\Delta > 0 $. Since $N >0 $ if and only if  $ q_0 > 1/\lambda $ the last remark  follows from Remark \ref{rema}. 
\fin 

\noindent We shall need the following easy arithmetic:
\begin{Lemm} \label{lemn}
 Assume $ \lambda > 9 $  then  $\sigma > \sqrt{\rho} $ and $ 1/\lambda < \sigma -\sqrt{\rho}$.
\end{Lemm}  
 \noindent{\it Proof}.- To prove the first inequality: 
$$
\lambda > 0  \Longleftrightarrow(\lambda -3)^2 >( \lambda -9)(\lambda -1) \Longleftrightarrow 
(\lambda -3) > \sqrt{ (\lambda -9)(\lambda -1)}. 
$$
 \noindent
To prove the second inequality: 
$$ 
 \begin{array}{cc} 
2 <  &  ( \lambda -3) - \sqrt{ ( \lambda -9)(\lambda -1) } \\
\Longleftrightarrow  
             &  (\lambda -5) >\sqrt{( \lambda -9)( \lambda -1)}   \\
 \Longleftrightarrow  &   ( \lambda -5)^2 > \lambda^2 -10 \lambda + 9 
\end{array} 
$$ 
\noindent and the 
last inequality is always true. 
\fin 
\begin{Cor}\label{lemo} Let $ q_0 < 1/ \lambda $. 
$ q_0 + q_2 <1 \Longleftrightarrow  -N   +  \Delta > 0  \, \, \mbox{or} -N - \Delta >0 $.
\end{Cor} 
\noindent {\it Proof}.- By using the proof of Remark \ref{rema} $ q_0 + q_2 < 1 \Longleftrightarrow 
-(\lambda q_0 -1)( q_0 +1 )   + \epsilon  \Delta  < 2( 1-\lambda q_0) \Longleftrightarrow   +  \epsilon \Delta <
-N \Longleftrightarrow  0 < -N \pm \Delta  $ with  $ \epsilon \in \{ \pm 1 \} $.  
\fin
\begin{Prop}  Assume $\lambda > 9 $ and $ \lambda q_0 <1 $ then
$ q_0 + q_2 <1 $.
 \end{Prop}
\noindent{ \it Proof}.-~Clearly $\lambda q_0 <1 \Longleftrightarrow  \lambda q_0 -1<0 $ thus $N<0 $. 
By lemma \ref{lemu}  to prove $ \omega > 0 $ we verify one of  the inequalities : $ q_0 < \sigma-
\sqrt{\rho} $ which is true by lemma \ref{lemn} hence $ 1/\lambda < \sigma - \sqrt{\rho} $. 
By Prop.~\ref{imp} one of the inequalities is satisfied thus  $\Delta^2 > 0 $ and $ -N + \Delta  = -N + \sqrt{-N \lambda \omega}>0 $ which is one of the inequalities to be satisfied
in Corollary \ref{lemo}. 
\fin 
\begin{Prop}\label{cota} Assume $\lambda > 9 $ and $ q_0 \in ( \sigma - 
\sqrt{\rho}, \sigma + \sqrt{\rho}) $ then $ q_0 > 1/\lambda, \,\, q_0 + q_2 < 1 $. 
\end{Prop}
\noindent{\it Proof}.-~If $\lambda >9$  then $ 1/\lambda < \sigma- \sqrt{\rho} $ by 
lemma \ref{lemn} but then  $ q_0 > 1/\lambda $. In this case $ N > 0 $ hence 
by lemma \ref{lemp}  the inequality  $ q_0 + q_2 < 1$ is true.  
\fin
\begin{Obs}
The cases $\lambda =1, 9 $  in the above proposition give  no criteria to find $ q_0 $ as a solution to equation \ref{ecuaci}. This can be explained as follows: fix  the  affine plane $H= \{ g= q_0 + q_2 + q_4 -1 = 0 \} $ in the $ \R^3 $
defined by  the coordinates $ q_0, q_2, q_4 $ and define for each $ \lambda \in  \R_{ \ge 0 }$:
$$
f_{\lambda}=  q_2 q_4 + q_0 q_4 + q_0 q_2 - \lambda q_0 q_2 q_4
$$
 which is   a surface in the ${\Af}_3 $ defined by the coordinates $ q_0, q_2, q_4 $. Let $ C_{ \lambda} =
H \cap \{ f_{\lambda} = 0 \} $.
\end{Obs}
\begin{Prop} The linear system of curves 
$\{ C_{\lambda}= \{ f_{\lambda} = g = 0 \} \}_{\lambda \in \R_{ \ge 0}} $ is always smooth 
except for    $ \lambda =1, 9 $ with singularities:
$$
\mbox{ Sing}( C_9 ) = \{ ( 1/ 3, 1/3, 1/3 ) \}, \,\, \mbox{ Sing} ( C_1) = \{ Q=  ( -1,1,1) \}.
$$
\end{Prop}  
\noindent{\it Proof}.- To simplify the computations let  $ \partial_i f = \partial f / \partial q_i $. Recall that
$ f_{ \lambda} = q_2 q_4 + q_0 q_2 + q_0 q_4 - \lambda q_0 q_2 q_4 $ and $ -g= 1 - q_0 -q_2 -q_4 $. Note that  $ \partial_i f
= q_j + q_k - \lambda q_j q_k  $ for  $ i, j, k \in \{ 0,2,4 \} $ with $ i \neq j \neq k $ and $ \partial_i g = -1 $. The jacobian matrix  is 
$$
\left(   \begin{array}{ccc}
 \partial_2 f_{\lambda}              &   \partial_0 f_{\lambda}       &       \partial_4  f_{\lambda}      \\
            -1           &    -1         &    -1
\end{array}
\right).
$$
It is of rank less than one if and only if:
 $$
0 = (\partial_i - \partial_j) f_{\lambda}  = 0 \, \,  i,j \in \{ 0,2,4 \}.
$$
 It follows  that : $ ( q_i -q_j ) ( 1 -\lambda  q_k ) = 0 \,\, \mbox{for some} \,\, \, i,j,k \in \{ 0,2,4 \} $. Due to the symmetries  of the indices
 it is enough to consider the following cases:

 \noindent{\bf I}. $ 1- \lambda q_4 = q_0 - q_4 = q_2 -q_4 =0 $  .  $ q_2 = q_4 = q_0 $.  $ q_4 = 1/ \lambda $. From  $ 0 =  g  =
 1 - 3 q_0 $. Therefore  $ q_0 = 1/3 $. Substituting  $ 0 =  \fl  =  3 - \lambda / 3 $. Thus $ \lambda = 9 $.

\noindent{\bf II}. $ 1 -\lambda q_4 = 1 -\lambda q_2  = q_2 -q_4 = 0 $.~Therefore $ q_2 = q_4 = 1/ \lambda $.
 $ g  = 1-2/ \lambda - q_0 = 0 $.
Substituting  $ q_0 = 1 - 2/ \lambda $ in  $ \fl  = 0 $ one obtains $ 0 = \fl =  
(1 + \lambda  -2)/{\lambda} = 0 $. Thus
 $ \lambda = 1$.

\noindent{\bf III}.$ 1-\lambda q_4 = 1-\lambda q_2 = 1-\lambda q_0  = 0 $.~Therefore $ q_0 = q_2 = q_4 =  1 / \lambda $.~Hence
$ \fl  =  2 / \lambda^2  = 0 $ which is impossible.

\noindent{\bf IV}.  $ 0 =  q_2 -q_0 =  q_0 - q_4 =  q_2 - q_4 $.~It follows that  $ q_0 = q_4  = q_2 $.~Hence $ g  =  1 - 3 q_0 = 0 $.
Substituting $ q_0 = 1/3 $  in $\fl $ one obtains  $ \fl  =  3/9 - \lambda / 27 = 0$ therefore $ \lambda = 9 $.
\fin

\noindent
Fix once again the $\R^3$ defined by the coordinates $ x,y,z$.
\begin{Lemm} 
The equation of $ \{C_{ \lambda} = \{ 0= \fl(x,y,z) = 1 -( x + y +z) \} \}_{\lambda \in {\R}_{>0}} $ can be written for $ \lambda  = 9 $ as :
$$
\begin{array}{cc}
f_9  &  =   9 ( x^2y + xy^2) + x + y -10 xy - (x^2 + y^2), \\
f_1  & =   x^2y + xy^2 + x + y  -2 x y - ( x^2 + y^2).
\end{array}
$$
\end{Lemm}
The proof of the lemma is left as an easy exercise.~For the next lemma recall that 
 $ P= ( 1/3, 1/3, 1/3) $ is a singular point for $C_9 $ and $ Q= \mbox{Sing}(C_1)$.
\begin{Lemm} The equation for the tangent cone of $C_9 $ (resp.~of $C_1$)   passing through  $P$(resp. $ Q$ of $C_1 $) is
$ ( x- 1/3)^2 + ( x-1/3)(y-1/3) + ( y-1/3)^2 $ (resp. $ 4(y-1)(x + y) $).
In particular, $P$  (resp. $Q$) is a non-ordinary double point of $C_9 $ at $P$ (resp. $Q$ of $C_1$).
\end{Lemm}
\noindent{ \it Proof}.- The second partial derivatives  at $P$ are  given as
 $ \partial^2_x f_9 = -2 + 18 y $, $\partial^2_{xy} f_9 = -10 + 18 x + 18 y $,
 $ \partial^2_{y} f_9 = -2 + 18x $. Summarizing : $ \partial^2_x f_9(P) = 4, \partial^2_y f_9(P) = 4,
\partial^2_{xy}f_9(P)= 2$. The equation for the tangent cone at $P$ is then
$$
4( x-1/3)^2 + 4(x-1/3)(y-1/3) + 4(y-1/3)^2.
$$
The calculation for $C_1$ can be done analogously.
\fin

\noindent
 \begin{Obs}A direct computation shows that  our $ f_9 =0 $ is irreducible over $\R$ .~Under the linear change of coordinates $ u= x-1/3, v= y-1/3 $ the equation for $ f_9 = 0 $ is transformed to
$ f_9 = 9 uv( u+ v) + 2 uv + 2( u^2 +v^2) $.~Under this linear change of coordinates the cubic  curve $C_9 $ is transformed to
a real cubic  with isolated singularity at the origin which is to be expected from the
  classification of irreducible cubic curves over the real number field.
\end{Obs}
\section{ The  32 lines  on the quartic surface.}\label{superf}
  We start our dicussion with some  well known  facts on Line Geometry and Group theory.~Fix the three dimensional real projective space $\R {\Pr}_3 $ with coordinates $ z_0, z_1, z_2, z_3$.~Introduce  the quartic  surface   $ X_f = \{ f(z_0, z_1, z_2, z_3) = 0 \}$ given  
by the homogeneous polynomial $f$ of degree four in the variables $ z_0, z_1, z_2, z_3 $.~A line in $ \R {\Pr}_3 $is generated by a two plane in $\R_4$ represented
 by  a two by four matrix
 $$
\left( 
\begin{array}{cccc}
 1 & 0 & * & *   \\
 0 & 1 & * & *
\end{array}
\right)
 $$
(There are in fact 6 ways  of choosing the matrix with  this property !)
In studying problems related to the study of  the  Geometry of hypersurfaces of the grassmanian of lines  in $\Pr_3$, known as 
Line Geometry, there exist various choices for the coordinates to use.
In studying the line complex it was typical in the ninteenth and  the early twentieth century  to introduce special cooordinate systems such as  elliptic, pentaspherical
coordinates (c.f. [\cite[chap. VIII, \S 130 and chap. XII, \S 221]{Je}.~A coordinatefree  approach uses  only linear
algebra to characterize properties of the line complex (c.f. \cite[chap.~6]{GH}). Here we choose the first approach due to
 the symmetries of our problem ( c.f.\cite{BN} for  the relation of this approach to polarized abelian surfaces).~More precisely let
$$
\Lambda =  \left(   
\begin{array}{cccc}
  z_0  & z_1 & z_2 & z_3       \\
 z_0^{'} & z_1^{'}  & z_2^{'}  & z_3^{'}
\end{array}
\right)
 $$
 be a two by four matrix and  introduce the following coordinates
 $$
p_{i,j} = z_i z_j^{'} - z_j z_i^{'}   \, \,  \mbox{i,j} \in \{ 0,1,2,3 \}
$$
where $ i \neq j$.~These are the Pl\"ucker coordinates or P-coordinates for short ;
In such coordinates the condition that a matrix $\Lambda$ defines a two-plane is given as:
$$
p_{01}p_{23} -p_{02}p_{13} + p_{03}p_{12} = 0.
 $$
\noindent In terms of multilinear algebra this is nothing else than giving a two form $\omega$ such that $ \omega \wedge \omega =0 $
(c.f.~\cite[chap.~1]{GH}).~In the coordinates $\{ p_{i,j} \} $ this is a hypersurface of degree two ( the Pl\"ucker quadric).~Let  $H^t $ be the subgroup of $\mbox{SL}(4, \R) $ spanned by the transformations
$$
 \begin{array}{cc}
 \sigma_1  = \left(  
             \begin{array}{cccc}
                0 & 0 & 1 & 0  \\
                0 & 0 & 0 & 1 \\
                1  & 0 & 0 & 0 \\
                  0 & 1 & 0 & 0
              \end{array}
                \right)
    &  \sigma_2 = \left(  \begin{array}{cccc}
               0 & 1 & 0 & 0  \\
              1 & 0 & 0 & 0  \\
                0 & 0 & 0 & 1 \\
               0 & 0 & 1 & 0
                 \end{array} 
\right) 
                     \\
\tau_1 = \left( 
            \begin{array}{cccc}
          1 &  0   &  0 &  0    \\
          0 & 1 & 0 & 0  \\
         0  & 0 & -1 & 0 \\
         0 & 0 & 0 & -1
          \end{array} 
              \right) &
\tau_2= \left( 
                      \begin{array}{cccc}
                        1 &  0  &  0 &  0  \\
                      0    &  -1 & 0 & 0  \\
                      0  & 0 &  1 &  0 \\
                      0 & 0 & 0 &   -1
                           \end{array} 
                              \right)
\end{array}
$$
 which  satisfy the relations:
$$
\sigma_i^2 = \tau_i^2 = \mbox{id},  \, \, \sigma_i \tau_i = - \tau_i \sigma_i
$$
for $ i=1,2 $.~One obtains a central exact sequence of groups:
$$
1 \rightarrow \{ \pm 1 \} \rightarrow H^t \rightarrow G' \rightarrow 0
$$
where $ G' \simeq {{\Z}_2}^4 $.~The {\it  explicite} action of $H^t $ on $f $ for   a  polynomial  
as above on the P-coordinates  is induced  by the  usual linear action induced on the polynomials of degree four,
 in particular
$$
 \sigma(z_0 z_1 z_2 z_3 ) = z_0 z_1 z_2 z_3  \, \, \mbox{for all} \, \, \sigma \in H^t .
$$
  The action of $ H^t$, the (unique up to a constant) Schroedinger
 representation of
degree four  on ${\R}^4$ induces  a representation on $ \wedge^2 \R^4 $  as  given  in \cite[p. 178]{BN}.  
\noindent Introduce the following coordinate transformation 
 in  $ {\R} {\Pr}_4 $ :
$$
\begin{array}{cccc}\label{matr}
x_0 & = p_{01} - p_{23} , \, \, x_2 & = p_{02} + p_{13}, \,  x_4 & = p_{03} -p_{12},  \\
 x_1 & = \sqrt{-1}( p_{01} + p_{23}), x_3 & = \sqrt{-1}(p_{02} - p_{13}), x_5 & = \sqrt{-1}(p_{03} + p_{12}).
\end{array}
$$
 Note that
\begin{equation} \label{pluc}
x_0^2 + x_1^2 + x_2^2 + x_3^2 + x_4^2 + x_5^2 =-2( p_{01}p_{23} - p_{02}p_{13} + p_{03}p_{12})
\end{equation}
which is the equation for the Pl\"ucker quadric, which parametrizes  the set of lines in $\R {\Pr}_3$.~The importance of the coordinates $\{ x_i \} $ called the Klein coordinates
 ,~hereafter K-coordinates,
  is that these are eigenfunctions for the action of the Heisenberg group $H^t$ on them.
 These eigenfunctions  are useful for studying  the configuration of lines
as in  \cite[section 4]{BN}.~Equation ~\ref{pluc}  parametrises the set of lines  in ${\R} {\Pr}_3 $.~By inverting   the transformation induced by equation~\ref{matr}  we obtain
 the coordinates for the P-coordinates  in terms of the K-coordinates as presented
in the introduction.~Using equation~\ref{matr}  the statement of the PROBLEM stated in section \ref{algor} is formulated equivalently as:

\noindent {\it The line with P-coordinates}
$ \{ p_{i,j} \} $
 {\it defines a real line  if and only if} $ x_0, x_2, x_4$ {\it are real and } $ x_1, x_3, x_5$ {\it purely imaginary}.

\noindent  The P-coordinates can be expressed  in terms of the $\{ q_i \} $-coordinates as: 
 \begin{eqnarray}
 \left( \begin{array}{c} 
p_{01} \\
p_{03} \\                                                          
 p_{13}
\end{array}    
\right)   
  &  =     
 \left( \begin{array}{c}
 \sqrt{q_0} + \sqrt{q_1}  \\ 
\sqrt{ q_4} + \sqrt{q_5} \\
 \sqrt{q_2} -\sqrt{q_3}
\end{array}
\right) 
, \\ \nonumber 
\left( \begin{array}{c}
p_{02}  \\   
    p_{23} \\        
 p_{12}       
\end{array}               
\right)
  & =  \left( \begin{array}{c}
  \sqrt{q_2} +  \sqrt{q_3}  \\
  -( \sqrt{q_0} - \sqrt{q_1} ) \\
  -( \sqrt{q_4} - \sqrt{q_5} )  
  \end{array}
\right). 
\label{recta}
\end{eqnarray}    
Let us consider   a line $ l \subset X$ with coordinates $\{ p_{i,j} \} $ such that   $p_{01} \neq 0 $ ( which is  equivalent to 
  $\sqrt{q_0} \neq -\sqrt{q_1} $).
  The  two points with coordinates 
 $$
\left( 
\begin{array}{c}
p_a  \\   
p_b
\end{array} 
\right) 
= \left( 
\begin{array}{cccc}
1 & 0 &  x &  y   \\
0 &   1 &   u &   v 
\end{array}
\right)
$$
in $ \R {\Pr}_3 $  are determined  by the P-coordinates as: 
$$
\begin{array}{cccc}
  u & = p_{02},    & \quad  v  & = p_{03},   \\
y & = -p_{13},  & \quad  x & = -p_{12}.
\end{array}
$$
The line spanned by $p_a, p_b$ will be  by definition $l_{a,b}$.

\noindent Let us consider  quartic surfaces which are invariant under the action
 of $H^t$.~We fix  once again  the  quartic surface $ X= \{ f = 0 \}$
then
$ \sigma f = \lambda_{\sigma} f $ for all $ \sigma \in H^t $ such that $ \lambda_{ \sigma} \in \{ \pm 1 \} $. 
One sees   inmediately (\cite[ Prop. 4.1.1 ii)]{Nie} )  that the set of real quartic forms
 invariant under $H^t$ is a real vector space of dimension five.~By fixing such a form it depends on the real values
 $\lambda_0, \lambda_1, \lambda_2, \lambda_3, \lambda_4$  which are
 the coefficients of that form.~If  a form of degree four contains a line  it imposes
 five linear conditions on the coefficients.~More precisely, the   restriction of   $f$ to a line $l_{a,b}$ induces five
 linear equations in the coordinates of the two points on the line and the coefficients 
   $\lambda_i$ for  $i=0, \ldots, 4 $.
 It
follows (\cite[\S 4.3]{BN})  that each line  $l_{ab}$  with K-coordinates satisfying  equation~\ref{proble}  belongs
 to  at least one quartic surface defined by a quartic form invariant under $H^t$.~Therefore the five linear equations
are linearly dependent and  projectively four of these equations determine the set
of coefficients $\{\lambda_i \} $ of the quartic form which  depend on
 the equivalence class of $ l_{a,b} $ modulo the action of $H^t$.~We  will assume  in section \ref{ejemp}  that 
 $\lambda_4 = 1$. 
\noindent The  Maple program given in section \ref{ejemp} allows us to determine 
two other points of the line given  two of the points of the line $l_{ab}$  hence by the previous discussion determine the coefficients $\lambda_i $ for $ i=0, \ldots, 4 $ of the quartic surface.~If an $H^t$ invariant quartic contains a line $l$ then
the orbit of $l$ under $G'$ contains 16  skew lines.~In the K-coordinates $\{ x_i \}$  the involution
$$
^{´}:(x_0: x_1: x_2: x_3: x_4: x_5) \mapsto    ( -1/ x_0: 1/x_1: 1/x_2: 1/x_3: 1/x_4: 1/x_5)
$$
which is well defined away from the fourfolds $\{ x_{i} =0 \}$ applied to $l$ gives a line $l'$. Writing the
P-coordinates associated to this line (for this let
$  q' =  -1/ \sqrt{q_0} -1/ \sqrt{q_1} $ ):  
\begin{eqnarray}
 \left( \begin{array}{c}
  p_{01}'       \\                                                        
 q' \cdot p_{03}'  \\                    
 q' \cdot  p_{13}'
\end{array}    
\right)   
  & =     
 \left( \begin{array}{c}
 1  \\ 
1/\sqrt{q_4} -1/\sqrt{q_5} \\
1/\sqrt{q_2} + 1/\sqrt{q_3},
\end{array}
\right) 
, \\ \nonumber 
\left( \begin{array}{c}
 q' \cdot p_{02}'  \\   
  q' \cdot   p_{23}' \\        
  q' \cdot p_{12}'       
\end{array}               
\right)
  & =  \left( 
\begin{array}{c}
  1/\sqrt{q_2}- 1/ \sqrt{ q_3} \\
    1/\sqrt{q_0} -1/\sqrt{q_1}\\
 -1/\sqrt{q_4} + 1/\sqrt{q_5}  
  \end{array}
\right). \label{rectat}
\end{eqnarray}    
The previous P-coordinates allows us in the  Maple program  
 given in section \ref{ejemp}  to determine very explicitely the 
parametric equation of the transerval lines once given  the line $l_{ab} $.
  Using the K-coordinates  
  $l'$ cuts  exactly 10 lines  of the Heisenberg orbit of $l$ (c.f. \cite[prop. 4.2 b)]{BN}). By reasons of degree, $l' \subset X$.
 Hence the orbit of  $l'$  under $G'$ is contained in $X$. The two orbits of lines that are
 in X can be grouped as  the ``even '' lines  as  those  having an even
 number of  minus   signs in its K-coordinates  and the ``odd'' lines  as having an
 odd number of minus    signs in its K-coordinates
 ( in fact  using these K-coordinates we  have studied in \cite[prop.~4.2]{BN}  grouptheoretical properties of the  configuration of  these lines).
  
\noindent Classically, the set of odd
 and even lines form a configuration of type $(16_{10}, 16_{10})$ on the quartic surface 
 forming  a {\it double } 16 ( as is referred in  \cite[{\it loc. cit.}]{BN}).  By  \cite{Se} a
 smooth quartic surface  can contain at most 48 lines
so if it contains at  least  32  lines
 it  must contain exactly 32 lines. \e
\section{A brief description of  the program }
\noindent The program written in Maple IV  release 4  defines  the following  local variables:
$$
R, S, q_1, q_0, M, N, Sq_1, Sq_2, q_2, q_3, q_4, q_5,M_1, N_1,  rr, ss, zz, com, ww, m, n, K, quarn
$$
and  the global variable $d$.~It uses lemma \ref{valor} given in section \ref{algorn}.~The initial value for the program is  $ \lambda $  ( in the program  this value is  $d $ ).~It then
evaluates  $ q_0  $ and $ ( R, S) $ using prop.~\ref{valor}, this calculation is performed in the procedures Var, Vas. In order to evaluate  the variable $ q_2 $  one needs to introduce  the variables
$M, N $  in terms of $d $  and  $ q_0 $  and finally evaluate  $ Sq_1 = \sqrt{ MN} $.~The algorithm then evaluates  
 the  positive root  of $ q_2 $ in terms of  $ N, Sq_1 $ and $  dq_0 - 1 $.~It follows that  $ q_4 = 1- ( q_0 + q_2) $.~Analogously,  if one gives  the value for $ q_1$  the program evaluates  $ q_3, q_5 $  in exactly the same way.

\noindent In order to  draw the lines one evaluates the parametric  equation  for the lines.~The procedure to evaluate them in the program
has been given in section \ref{superf}.~We introduce the  variables $ rr, ss $ in terms of  $ q_4, q_5, q_2, q_3 $.~The  last variables are
used to give the parametric equation of the line $ l $ as given  by equation~\ref{recta}.~The orbit of $l $ under $ H^t $ is evaluated using the procedure
 graf.~The parametric equation of the  transversal  to $ l $ is evaluated using the variables  $ m,n $.~For these  one uses the
variables $ rr, zz, ww $ and uses equation~\ref{rectat}.~The  coefficients for the  quartic surface  are evaluated by means of a matrix array  using the
elementary method of the previous section.~The equation of the quartic  surface in a fixed solid herafter named
the ``clipping solid''  is evaluated in a procedure and saved as the variable quarn.~The  procedure graf  used to draw the surface (resp.~
the lines lying on this surface) is   based on the command (resp.~spacecurve) display3d of  the library plots of Maple. 
\section{An Example}\label{ejemp}
\noindent
To illustrate the above theory  we will introduce a typical example.~The value for this example was chosen from a series of  43 test examples for $d \in [9.3, 18.00]$ on very different intervals, on rectangular grids of  $ 30^3= 27,000 $ points  (using the command  grid).~The program was run on a Samsung pentium II.~The average compilation time per image for 
these examples averaged  to three minutes.~We introduce the value for  $d$ which in this case is $18.00$ (this is the upper bound!).~The necessary condition for $ q_0$ is given in lemma~\ref{valor}.~For practical purposes we introduce two small routines written in Maple  to evaluate $q_0, q_1$ in terms of $d$.~For this example, $ q_0= 0.4168, q_1 = 0.1713 $.~Using the value for $q_0$ and two other procedures written on the program one obtains the values 
  $ M =2.124999680,  N = 3.792199680 $.~Using $ Sq_1 = 2.83874$ one evaluates   $ q_2 = 0.5101,  q_4 = 0.0731$. The values for
$M_1, N_1, Sq_2, q_3, q_5, rr,ss, zz , m , ww ,n  $ are (respectively)
  $1.041313580, 1.762513580, 1.3408,
0.7364,0.0923,-0.0334,1.5724,0.14393,\break -0.10242,0.57418,-0.059153$.~Using the elementary theory described in section \ref{superf}, one obtains  a four by four matrix  $A$  with entries polynomials of degree  four  in the variables  $x, y, u,v $  and a four by one vector 
such that the extended matrix $ M = ( A, b ) $  gives the equation $ M \cdot \lambda = 0 $ meaning that the quartic surface with vector coefficients  $ \lambda $ contains a line.~Written  in  non-homogeneous  form one obtains : $ A \cdot \lambda = b $.~The procedure mpoly  substitutes the  variables  $ x, y,u,v $ for the obtained values  $ rr,ss,zz, ww $.~For this example  $ \lambda = ( -0.366,-1.44, 0.614,0.526) $.~We have chosen the
clipping solid as $\{  x + v + y + u = 0 \} $. In this solid  the equation for the
quartic polynomial is given by:
\begin{eqnarray*}
 \lefteqn{quarn  = 0.320 x^4 + 0.496 y^4 + -3.612 u^4 } \\  
   & &       -2.796  (x^2 y^2 + x^2 u^2 + y^2 u^2 ) +   uy ( -u^2 7.224 + 0.992 y^2) \\
   & & + xy(  0.640 x^2   + y^2  0.992)   
   + xu( x^2 0.640  - u^2 7.224)    \\ 
 &   &    +xyu( -y 5.936 - 14.152 u -x 6.288).  
\end{eqnarray*}
\noindent In the image produced by the procedure  graf of the Maple program in all the examples, the line $l$ has been drawn in  red.~The remaining colors for the disjoint lines
are: blue, yellow, sienna, cyan, khaki, pink, turquoise, aquamarine, magenta, plum, violet, braun, green, navy, gold.~The ten transversals to $l$ have been drawn in color grey.~In the plot3d command of Maple one has to specify a pair of values $(u,v) $for the orientation i.e.~the direction from which the  object is to be viewed.~For  the  same group of test examples described in the interval mentioned at the beginning of this paragraph we tested  pairs $ u \in ( 0, 45) $ and $ v \in ( 0, 90)$ showing in each image at most 7 transversals, with no further improvement.~We illustrate an example at the end of  the references for the surface  defined by the above given polynomial  with a grid of 27,000 points allowing  us to see seven 
of the  transversal lines.~It is for $ d= 18.00 $ and  at $ (u,v)= (20,36)$. 

\begin{figure}[b]
\includegraphics[340,340]{sphere.eps}
\center
\caption{The quartic surface for $ d= 18.00$ at  $(20,36) $.}
\end{figure}

\noindent www.ifm.umich.mx/personal.html/Isidro Nieto 

\noindent{\tt inieto@zeus.ccu.umich.mx}

\noindent{\it Current address}:
 
\noindent I. Nieto 

\noindent Apartado postal 2-82 

\noindent 58040  

\noindent Morelia, Mexico.

\begin{thebibliography}{basura}
\bibitem [BN]{BN} Barth, W. and Nieto, I. Abelian surfaces of type (1,3) and quartic surfaces with 16 skew lines {\it J.~Alg.~Geom.},~{\bf 3}~(1994), pp.~173-222.
\bibitem [GH]{GH} Griffiths, P. and Harris, J.~{\it Principles of Algebraic Geometry},  N.Y. Wiley (1978).
\bibitem [Je]{Je} Jessop, C.M.{\it A Treatise on the Line Complex}, New York, Chelsea ( 1903).
\bibitem [Nie]{Nie} Nieto, I.~{\it Invariante Quartiken Unter der Heisenberg Gruppe T},~PhD Thesis, U. of Erlangen, (1989).
\bibitem [Nie1]{Nie1} Nieto, I. Examples of abelian  surfaces with polarization type (1,3)  in {\it Algebraic Geometry  and Singularities} ( ed. C. Lopez and N. Macarro), Progress in Mathematics, vol. 134, pp. 319-337, Birkhauser, Basel (1996).
\bibitem [Se]{Se} Segre, B. The maximum number of lines lying on a quartic surface~{\it Oxf. Quart. Journ.},~{\bf 14} (1943),~pp.~87-96.
\end{thebibliography}
\end{document}